\documentclass{article}

\usepackage{arxiv}

\usepackage{algorithm}
\usepackage{algpseudocode}

\usepackage[utf8]{inputenc} 
\usepackage[T1]{fontenc}    
\usepackage{hyperref}       
\usepackage{url}            
\usepackage{booktabs}       
\usepackage{amsfonts}       
\usepackage{nicefrac}       
\usepackage{microtype}      
\usepackage{graphicx}
\usepackage[numbers]{natbib}
\usepackage{doi}
\usepackage{amsmath,mathtools}
\usepackage[commentmarkup=todo]{changes}
\usepackage{todonotes}
\usepackage{rotating}
\usepackage{tikz}
\usepackage{pgfplots}

\usepackage{booktabs,multirow}
\usepackage{siunitx}
\definecolor{TolMutedBlue}{HTML}{332288}
\definecolor{TolMutedCyan}{HTML}{88CCEE}
\definecolor{TolMutedTeal}{HTML}{44AA99}
\definecolor{TolMutedGreen}{HTML}{117733}
\definecolor{TolMutedOlive}{HTML}{999933}
\definecolor{TolMutedSand}{HTML}{DDCC77}
\definecolor{TolMutedRose}{HTML}{CC6677}
\definecolor{TolMutedWine}{HTML}{882255}
\definecolor{TolMutedPurple}{HTML}{AA4499}

\hypersetup{
    colorlinks=true,
    citecolor=TolMutedBlue, linkcolor=TolMutedGreen,   urlcolor=TolMutedPurple,%
}

\definechangesauthor[name = {Mateusz Baran}, color = {TolMutedBlue}]{MB}
\definechangesauthor[name = {Ronny Bergmann}, color = {TolMutedRose}]{RB}
\definechangesauthor[name = {Patryk Przybysz}, color = {TolMutedGreen}]{PP}
\todostyle{RB}{inline, bordercolor=TolMutedRose, backgroundcolor=TolMutedRose!10}
\todostyle{MB}{inline, bordercolor=TolMutedBlue, backgroundcolor=TolMutedBlue!10}
\todostyle{PP}{inline, bordercolor=TolMutedGreen, backgroundcolor=TolMutedGreen!10}

\newcommand{\retr}{\ensuremath{\mathrm{retr}}}
\newcommand{\grad}{\ensuremath{\mathrm{grad}}}
\newcommand{\Hess}{\ensuremath{\mathrm{Hess}}}

\title{A Riemannian quasi-Newton algorithm for optimization with Euclidean bounds}

\author{Mateusz Baran\textsuperscript{\href{https://orcid.org/0000-0000-0000-0000}{\includegraphics[scale=0.06]{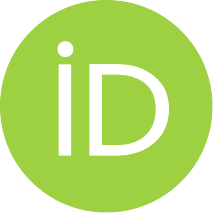}}}\\
	AGH University of Krakow\\
    30 Mickiewicz Ave., 30-059 \\
	Kraków, Poland \\
	\texttt{mbaran@agh.edu.pl} \\
	\And
	Ronny Bergmann\textsuperscript{\href{https://orcid.org/0000-0000-0000-0000}{\includegraphics[scale=0.06]{orcid.pdf}}}\\
    Department of Mathematical Sciences, \\
	Norwegian University of Science and Technology, \\
    NO-7041 Trondheim, Norway \\
    \\
    AGH University of Krakow\\
    30 Mickiewicz Ave., 30-059 \\
	Kraków, Poland \\
	\texttt{ronny.bergmannn@ntnu.no} \\
    \And
    Patryk Przybysz \\
    AGH University of Krakow\\
    30 Mickiewicz Ave., 30-059 \\
    Kraków, Poland \\
}

\renewcommand{\shorttitle}{Riemannian L-BFGS-B}

\hypersetup{
pdftitle={A Riemannian L-BFGS-B},
pdfauthor={Mateusz Baran, Ronny Bergmann, Patryk Przybysz},
}

\begin{document}
\maketitle

\begin{abstract}
	We propose a Riemannian limited-memory BFGS method for optimization problems with Euclidean bounds.
	The method combines a limited-memory quasi-Newton update in the tangent space with a Riemannian adaptation of the generalized Cauchy point strategy from classical L-BFGS-B, enabling efficient handling of Euclidean bounds while exploiting the geometric structure of the optimization domain.
	This setting is important in several applications, including covariance matrix estimation with bounded variance, neuroimaging, EEG signal classification, and other signal processing or computer-vision tasks that couple manifold variables with constrained Euclidean parameters.

	We provide a generic algorithmic framework and an implementation of the algorithm in the Manopt.jl library.
	Numerical experiments on benchmark problems indicate only minor reduction in performance on Euclidean problems compared to the classical L-BFGS-B method, while outperforming interior-point methods.
	Furthermore, the algorithm was tested on two mixed manifold and bounded Euclidean problems: amplitude-limited blind source separation with Gaussianity penalization and bounded-variance maximum likelihood common principal components analysis.
	The proposed method outperforms existing methods by several orders of magnitude.
\end{abstract}


\section{Introduction}
\label{sec:introduction}

Riemannian optimization is an important tool for solving optimization problems.
Many real-world applications, such as machine learning~\cite{FeiLiuJiaLiWeiChen:2025}, signal processing~\cite{ClosasOrtegaLesoupleDjuric:2024}, and computer vision~\cite{TuragaVeeraraghavanChellappa:2008,MettesGhadimiKeller-ResselGuYeung:2024}, involve optimization over manifolds.
For example, optimization over the manifold of positive definite matrices is common in gene expression data analysis~\cite{ChenKangJianYishiYunpengDonald:2018,KramerSchaferBoulesteix:2009}, neuroimaging~\cite{JuKoblerCollasKawanabeGuan:2026,BergmannTenbrinck:2018} and EEG signal classification~\cite{TibermacineTibermacineZouaiRabehi:2024}.
Grassmann manifold has also been used for simulations in nuclear physics~\cite{GongZhangYuanZhangXiong:2026}.
Geometric methods have also been more broadly applied in signal processing, for example in clustering~\cite{LabsirLesoupleTourneret:2026}, blind source separation~\cite{BouchardMalickCongedo:2018}, needle sensors~\cite{LezcanoIordachitaKim:2022} and microphone array signal processing~\cite{BarTalmon:2024}.

Adding bounds constraints to the problem allows for more flexible modelling of such problems~\cite{LiuWangZhao:2014}.
Constrains also offer a way to regularize a parameter estimation problem~\cite{FrankeHefenbrockKoehlerHutter:2024}.
In this work we introduce an efficient, generic algorithm for solving such problems to eliminate the need to develop purpose-built solvers for specific problems.

Quasi-Newton methods that effectively support bounds constraints, despite being known for a few decades, are still relevant in applications.
For example the L-BFGS-B algorithm has been recently shown to be effective in network traffic forecasting~\cite{LiuLeu:2025}.
Variants of L-BFGS-B can also be used as sub-solvers in the developing area of constraint learning~\cite{FajemisinMaragnoDenHertog:2024}.
On the other hand, handling bounds constraints in the Riemannian setting is currently limited to a few algorithms like augmented Lagrangian method (ALM)~\cite{LiuBoumal:2020}, exact penalty method (EPN)~\cite{LiuBoumal:2020} and interior point Newton method (IPN)~\cite{LaiYoshise:2024} which are generic but at the price of being potentially slower than algorithms that that are specialized to bounds constraints.

In this work, we consider problems of the form
\begin{equation}
    \min_{p=(p_{\mathrm{D}}, p_{\mathrm{M}}) \in D \times \mathcal{M}} f(x)
\end{equation}
where $D$ is the hypercube $[l_1, u_1] \times [l_2, u_2] \times \cdots \times [l_n, u_n]$ for some $n \geq 0$, $\mathcal{M}$ is a Riemannian manifold, $f\colon D \times \mathcal{M} \to \mathbb{R}$ is a smooth objective function, and $l_i, u_i$ are the lower and upper bounds for the variables $p_{\mathrm{D},i}$, $i=1,\ldots,n$.
We allow bounds to be infinite, that is, $l_i = -\infty$ or $u_i = \infty$ for some $i$.
Such problems often arise in statistics and machine learning, for example when estimating covariance matrix in a statistical model together with some other bounded parameters~\cite{DaVeigaMarrel:2012,SwilerGulianFrankelSaftaJakeman:2020,PensoneaultYangXiu:2020}, or covariance matrices with bounded covariances~\cite{ChenKangJianYishiYunpengDonald:2018}.

To solve this problem, we propose an approach that combines the Riemannian limited-memory BFGS algorithm~\cite{HuangGallivanAbsil:2015,HuangAbsilGallivan:2018} with our novel adaptation to the Riemannian setting of the generalized Cauchy point algorithm from the L-BFGS-B method~\cite{ByrdLuNocedalZhu:1995,ZhuByrdLuNocedal:1997,ByrdNocedalSchnabel:1994}.
We further provide an implementation within the framework of the Julia package Manopt.jl~\cite{Bergmann:2022:1}.
This framework allows to use the algorithm for arbitrary manifolds implemented within the interface defined in ManifoldsBase.jl,
for example all manifolds provided in the Manifolds.jl~\cite{AxenBaranBergmannRzecki:2023} library.

The remainder of this paper is organised as follows.
We first repeat the necessary terms and notation in Section~\ref{sec:preliminaries}.
Section~\ref{sec:methods} describes the proposed algorithm in detail.
Next, in Section~\ref{sec:results}, we present numerical experiments that demonstrate the effectiveness of our method on a set of benchmark problems.
Finally, Section~\ref{sec:conclusions} presents conclusions and discusses potential directions for future research.

\section{Preliminaries \& Notation}%
\label{sec:preliminaries}
In this section we introduce the main notation used throughout the paper.
We follow the textbooks~\cite{AbsilMahonySepulchre:2008,Boumal:2023} and refer to these for further details. We denote a Riemannian manifold by $\mathcal M$
and assume in the following that it is second-countable, Hausdorff, and complete.
At every point $p\in \mathcal M$ we denote the tangent space by $T_p\mathcal M$.
Elements of the tangent space are called tangent vectors and denoted by $X_p, Y_p$, where we omit the $p$ if it is clear from context.
The disjoint union of all tangent spaces it the tangent bundle $T\mathcal M$,
which is itself again a Riemannian manifold.
The Riemannian metric, is a family of inner products $\langle\cdot,\cdot\rangle_p$, $p \in \mathcal M$ on the tangent spaces that smoothly varies in $p$ and we denote the Levi-Cevita connection by $\nabla$.

The metric allows to introduce geodesics, i.\,e.~acceleration free curves
as well as the exponential map $\exp_p\colon T_p\mathcal M \to \mathcal M$,
which locally around $0\in T_p\mathcal M$ can be inverted and gives rise to the logarithmic map $\log_p$. Similarly the parallel transport $\mathcal P_{p,X}\colon T_p\mathcal M \to T_q\mathcal M$, $q = \exp_p(X)$ allows to “move” a tangent vector from one tangent space to another.
Numerically, these three operations might not be given in closed form for a certain manifold but require to solve an ODE. Then one can alternatively use
a retraction $\retr_p$ instead of an exponential map, a corresponding inverse retraction $\retr_p^{-1}$ instead of the logarithmic map, and a vector transport $\mathcal T_{p,X}$ to $q = \retr_p(X)$ instead of the parallel transport.

For a smooth function $f\colon \mathcal M \to \mathbb R$ the Riemannian metric
introduces the notion of a Riemannian gradient as its Riesz representer, that is,
via the unique tangent vector $\grad f(p) \in T_p\mathcal M$ such that $Df(p)[X] = \langle \grad f(p), X$ and the covariant derivative of the gradient introduces the notion of a Hessian: we have $\Hess f(p)[X] = \nabla_X \grad f(p)$

While in general, the product of two manifolds $\mathcal M_1 \times \mathcal M_2$ is again a Riemannian manifold using the usual product metric, note that the hypercube
 $D \coloneqq [l_1, u_1] \times [l_2, u_2] \times \cdots \times [l_n, u_n]$
 is not a Riemannian manifold.
 But it can be considered as a manifold with~\cite{Joyce:2010}. Then the product $D\times \mathcal M$ with a Riemannian manifold $\mathcal M$ also yields a new manifold with corners.
 For these, we introduce a projection operator $\operatorname{Proj}_{T_p (D\times \mathcal{M})} \colon T_p (D\times \mathcal{M}) \to T_p (D\times \mathcal{M})$ as $\operatorname{Proj}_{T_{(p_{\mathrm{D}}, p_{\mathrm{M}})} (D\times \mathcal{M})} (X_{\mathrm{D}}, X_{\mathrm{M}}) = (Y_{\mathrm{D}}, X_{\mathrm{M}})$ where
\begin{equation}
	\label{eq:projection_Y_D}
	Y_{\mathrm{D},i} = \begin{cases}
		0 & \text{ when } p_{\mathrm{D},i} = l_i \text{ and } X_{\mathrm{D},i} < 0 \\
		0 & \text{ when } p_{\mathrm{D},i} = u_i \text{ and } X_{\mathrm{D},i} > 0 \\
		X_{\mathrm{D},i} & \text{ otherwise},
	\end{cases}
\end{equation}
for each $i = 1, 2, \dots, n$.
We further introduce a Riemannian metric thereon, namely also the product metric $\langle (X_1, X_2), (Y_2, Y_2) \rangle_{(p_{\mathrm{D}}, p_{\mathrm{M}})} = X_1^{\top}Y_1 + \langle X_2, Y_2 \rangle_{p_{\mathrm{M}}}$ on $T_p (D\times \mathcal{M})$. Similarly, retractions, inverse retractions and vector transports
on this product can be considered component-wise.

\section{Methods}%
\label{sec:methods}

The general outline of the algorithm is similar to the Riemannian L-BFGS method~\cite{HuangAbsilGallivan:2018}
starting with an initial point $p_0$ on some manifold.

Then the idea of a single $k$th step is to solve the Newton equation as in \cite[Algorithm~5]{AbsilMahonySepulchre:2008}.
This is relaxed, since both setting up that linear system $\Hess f(p_k)[X_k] = -\grad f(p_k)$ as well as solving it might be expensive. Instead, one aims to find approximations $H_k \approx \Hess f(p_k)$ either
as a linear operator or as a coefficient matrix with respect to a certain basis of the tangent space at $p_k$.
Instead of approximating the Hessian, a popular alternative is to approximate the inverse $B_k \approx (\Hess f(p))^{-1}$. The advantage is, that computing the new direction $X_k$ then simplifies to applying the linear operator $B_k$ to the negative gradient. One prominent scheme here is the limited memory (L) variant of Broyden, Fletcher, Goldfarb, and Shanno (BFGS) which is based on differences of a certain limited amount of both iterates and gradients.
On manifolds, this turns into inverse retractions of iterates and differences of two gradients, where on of them is transported to the tangent space of the other.
Once a new direction has been found, the new iterate is obtained using a retraction as $p_{k+1} = \retr_{p_k}(\alpha_k X_k)$, where $\alpha_k$ is a step size, either as a fixed sequence of even employing some backtracking procedure.

In our setting, we begin with an initial point $p_0 \in D \times \mathcal{M}$.
At each iterate $p_k$, $k = 0, 1, \dots, k_{\mathrm{max}}$ a quadratic surrogate $m_k \colon T_{p_k} (D \times \mathcal{M}) \to \mathbb{R}$ of the objective is constructed:
\begin{equation}
	m_k(X) = f(p_k) + \langle \operatorname{grad} f(p_k), X \rangle_{p_k} + \frac{1}{2} \langle X, H_k[X] \rangle_{p_k},
\end{equation}
where $H_k \colon T_{p_k} (D \times \mathcal{M}) \to T_{p_k} (D \times \mathcal{M})$ is a positive definite approximation of the Riemannian Hessian of $f$ at $p_k$.
Section~\ref{sec:hessian_approximation} describes the details of construction of the Hessian approximation $H_k$ and its inverse using a limited memory approach, adapting the Euclidean approach to the Riemannian setting.

We compute the initial search direction by solving the subproblem $d_k = -B_k \operatorname{grad} f(p_k)$, where $B_k$ is the inverse of $H_k$.
It differs from the Euclidean L-BFGS-B method, where $d_k = (d_{\mathrm{D},k}, d_{\mathrm{M},k})$ is selected as $-\operatorname{grad} f(p_k)$.
After computing the initial search direction, we solve the subproblem of finding the first minimum of $m_k$ along the piecewise linear function $d_{\mathrm{PL}}(t) \colon [0, \infty) \to T_{p_k} (D \times \mathcal{M})$, $d_{\mathrm{PL}}(t) = (d_{\mathrm{PL,D}}(t), t d_{\mathrm{M}, k})$ where $d_{\mathrm{PL,D}}$ is defined as
\begin{equation}
	d_{\mathrm{PL,D}}(t)_i = \begin{cases}
		l_i - p_{\mathrm{D},i} & \text{ if } p_{\mathrm{D},i,t} < l_i \\
		t d_{\mathrm{D},i} & \text{ if } l_i \leq p_{\mathrm{D},i,t} \leq u_i \\
		u_i - p_{\mathrm{D},i} & \text{ if } p_{\mathrm{D},i,t} > u_i
	\end{cases}
\end{equation}
where $p_{\mathrm{D},i,t} = p_{\mathrm{D},i} + t d_{\mathrm{D},k,i}$.
Let the first minimizer of $q_k(t) = m_k(d_{\mathrm{PL}}(t))$ be attained at $t_{*,k}$.
The generalized Cauchy direction is then defined as $d_{\mathrm{GCD},k} = d_{\mathrm{PL}}(t_{*,k})$.
To find the minimizer we adapt the generalized Cauchy point algorithm from~\cite{ByrdLuNocedalZhu:1995} to the Riemannian setting; see Section~\ref{sec:gcd} for details.
Furthermore, the procedure returns the maximum allowed stepsize $t_{\mathrm{GCD},\max,k}$.
Note that in the Euclidean L-BFGS-B method, $\exp_{p_k}(d_{\mathrm{GCD},k})$ is called the generalized Cauchy point~\cite{ByrdLuNocedalZhu:1995}.

Next, we perform a line search along $d_{\mathrm{GCD},k}$ to ensure sufficient decrease in the objective function.
The line search operates within the $[0, t_{\mathrm{GCD},\max,k}]$ interval to find a suitable step size.

\subsection{Hessian approximation}%
\label{sec:hessian_approximation}

In this work we construct limited memory approximations to both the Hessian $H_k$ and its inverse $B_k$ simultaneously. This extends the standard Riemannian L-BFGS approach from~\cite{HuangGallivanAbsil:2015,HuangAbsilGallivan:2018}.

At each iteration $k=0, 1, \dots, k_{\mathrm{max}}$ we define
\begin{equation}
	s_k = \mathcal T_{p_k, \alpha d_k} (\alpha d_k) \quad \text{and} \quad y_k = \beta_k^{-1} \operatorname{grad} f(p_{k+1}) - \mathcal T_{p_k, \alpha d_k} (\operatorname{grad} f(p_k)),
\end{equation}
where $\beta_k$ is a scaling factor as defined in~\cite{HuangAbsilGallivan:2018}.

At iteration $k$ we have $\mu_k \leq \mu$ pairs of tangent vectors $(s_{k,i}, y_{k,i})$, $i = 1, 2, \dots, \mu_k$, where $s_{k,i}, y_{k,i} \in T_{p_k} D \times \mathcal{M}$ for $i = 1, 2, \dots, \mu_k$.
Initially $\mu_0 = 0$.
At iteration $k$, $k>0$, pairs $(s_{k,i}, y_{k,i})$ for $i \in {1, 2, \dots, \mu_k-1}$ are a subset of $(\hat{s}_{k,i}, \hat{y}_{k,i}) := (\mathcal T_{p_{k-1}, \alpha d_{k-1}} s_{k-1,i}, \mathcal T_{p_{k-1}, \alpha d_{k-1}} y_{k-1,i})$, $i = 1, 2, \dots, \mu_{k-1}$.
Elements such that $\langle \hat{s}_{k,i}, \hat{y}_{k,i} \rangle_{p_{k}} < \epsilon \lVert \hat{y}_{k,i} \rVert_{p_k}^2$ for some small positive constant $\epsilon$ are discarded to ensure the positive definiteness of the Hessian approximation~\cite{ByrdLuNocedalZhu:1995}.
If $\mu_{k-1} = \mu$ and no pair was removed by enforcing the positive definiteness condition, the oldest pair $(\hat{s}_{k,1}, \hat{y}_{k,1})$ is discarded.
Finally, the new pair is added at the end, $(s_{k,\mu_k}, y_{k,\mu_k}) := (s_k, y_k)$.

We further define $\rho_i = \langle s_{k,i}, y_{k,i} \rangle_{p_k}^{-1}$ and set the default Hessian scaling factor to $\theta_k = \langle y_{k,\mu_k}, y_{k,\mu_k} \rangle_{p_k}^2 \rho_{\mu_k}$.

To define the Hessian approximation, we first introduce some auxiliary functions.
Linear functions $W_{\mathrm{y},k} \colon T_{p_k} (D \times \mathcal{M}) \to \mathbb{R}^{\mu_k}$ and $W_{\mathrm{s},k} \colon T_{p_k} (D \times \mathcal{M}) \to \mathbb{R}^{\mu_k}$ compute coefficients of an input vector $X \in T_{p_k} (D \times \mathcal{M})$ with respect to the stored $y_{k,i}$ and $s_{k,i}$ vectors, respectively:
\begin{equation}
	W_{\mathrm{y},k}[X] = \begin{bmatrix}
		\langle y_{k,1}, X \rangle_{p_k} \\
		\langle y_{k,2}, X \rangle_{p_k} \\
		\vdots \\
		\langle y_{k,\mu_k}, X \rangle_{p_k}
	\end{bmatrix},
	\qquad  W_{\mathrm{s},k}[X] = \theta_k \begin{bmatrix}
		\langle s_{k,1}, X \rangle_{p_k} \\
		\langle s_{k,2}, X \rangle_{p_k} \\
		\vdots \\
		\langle s_{k,\mu_k}, X \rangle_{p_k}
	\end{bmatrix}.
\end{equation}
Next, these coefficients are used in a quadratic form
\begin{equation}
	\label{eq:H_W_k}
	H_{\mathrm{W},k}(\xi, c_{\mathrm{y},X}, c_{\mathrm{s},X}, c_{\mathrm{y},Y}, c_{\mathrm{s},Y}) = \theta_k \xi - \begin{bmatrix}
		c_{\mathrm{y},X} \\
		c_{\mathrm{s},X}
	\end{bmatrix}^{\top}
	M_k
	\begin{bmatrix}
		c_{\mathrm{y},Y} \\
		c_{\mathrm{s},Y}
	\end{bmatrix},
\end{equation}
where $M_k$ is the block matrix
\begin{equation}
\label{eq:M_k}
	M_k =
\begin{bmatrix}
		-D_k & L_k^{\top} \\
		L_k & Q_k
	\end{bmatrix}^{-1}
\end{equation}
composed from the following blocks:
\begin{equation}
	D_k = \operatorname{diag}(\rho_1, \dots, \rho_{\mu_k})^{-1},
\end{equation}

\begin{equation}
	Q_k = \theta_k \begin{bmatrix}
		\langle s_{k,1}, s_{k,1} \rangle_{p_k} & \cdots & \langle s_{k,1}, s_{k,\mu_k} \rangle_{p_k} \\
		\vdots & \ddots & \vdots \\
		\langle s_{k,\mu_k}, s_{k,1} \rangle_{p_k} & \cdots & \langle s_{k,\mu_k}, s_{k,\mu_k} \rangle_{p_k}
	\end{bmatrix},
\end{equation}
and
\begin{equation}
	L_k = \begin{bmatrix}
		0 & 0 & \cdots & 0 \\
		\langle s_{k,2}, y_{k,1} \rangle_{p_k} & 0 & \cdots & 0 \\
		\vdots & \vdots & \ddots & \vdots \\
		\langle s_{k,\mu_k}, y_{k,1} \rangle_{p_k} & \langle s_{k,\mu_k}, y_{k,2} \rangle_{p_k} & \cdots & 0
	\end{bmatrix}.
\end{equation}

Note that the invertion in Eq.~\eqref{eq:M_k} can be performed efficiently using the following formula:
\begin{equation}
	M_k = \begin{bmatrix}
		-D_k^{-1} - D_k^{-1} L_k Q_k^{-1} L_k^{\top} D_k^{-1} & D_k^{-1} L_k Q_k^{-1} \\
		Q_k^{-1} L_k^{\top} D_k^{-1} & Q_k^{-1}
	\end{bmatrix},
\end{equation}
where only one inversion of a nondiagonal $\mu_k \times \mu_k$ matrix $Q_k$ is required.

Finally, Hessian value $\langle X, H_k [Y] \rangle_{p_k}$, where $H_k \colon T_{p_k} D \times \mathcal{M} \to T_{p_k} D \times \mathcal{M}$ is then defined as
\begin{equation}
	\langle X, H_k[Y] \rangle_{p_k} = H_{\mathrm{W},k}(\langle X, Y \rangle_{p_k}, W_{\mathrm{y},k}[X], W_{\mathrm{s},k}[X], W_{\mathrm{y},k}[Y], W_{\mathrm{s},k}[Y]).
\end{equation}

For the generalized Cauchy direction algorithm we only need to compute the quantities of the form $\langle X, H_k[X] \rangle_{p_k}$, $\langle e_b, H_k [e_b] \rangle_{p_k}$, $\langle e_b, H_k [X] \rangle_{p_k}$ as well as the function $H_{\mathrm{W},k}$, where $e_b$ for $b = 1, 2, \dots, n$ is the tangent vector at $p_k$ that has $b$th standard basis vector as the first component and the zero vector as the second component.

\subsection{Generalized Cauchy direction}
\label{sec:gcd}

The input of the generalized Cauchy direction (GCD) algorithm is the current iterate $p$, the gradient $\operatorname{grad} f(p) = g = (g_{\mathrm{D}}, g_{\mathrm{M}})$, the search direction $d = -B \operatorname{grad} f(p) = (d_{\mathrm{D}}, d_{\mathrm{M}})$, and the Hessian approximation $H$.
We drop the outer iteration index $k$ for clarity.
Compared to the original generalized Cauchy point algorithm from~\cite{ByrdLuNocedalZhu:1995}, the output is the direction of descent $d_{\mathrm{GCD}}$, a status flag indicating whether the GCD was found limited by the bounds (and so the subsequent line search needs to operate within the $[0, t_{\mathrm{GCD},\mathrm{max}}]$ interval, where the upper bound is also computed), unlimited (line search may opt to perform bracketing) or not found (in which case the Hessian approximation needs to be discarded).

We split the logic into a part that is independent of the Hessian approximation and a part that depends on $H$ without increasing computational complexity.
This makes it easier to adapt the algorithm to different Hessian approximations.
Additionally, the algorithm works generically for any descent direction $d$, not only for the negative gradient direction.
Despite these extensions, the computational complexity remains the same as in the original generalized Cauchy point algorithm.

The first step of the algorithm is to compute the sequence of breakpoints $t_i$.
For each $i = 1, 2, \dots, n$ we compute
\begin{equation}
	t_i = \begin{cases}
		\frac{p_{\mathrm{D},i} - l_i}{d_{\mathrm{D},i}} & \text{ when } d_{\mathrm{D},i}
		 < 0 \\
		\frac{p_{\mathrm{D},i} - u_i}{d_{\mathrm{D},i}} & \text{ when } d_{\mathrm{D},i} > 0 \\
		\infty & \text{ when } d_{\mathrm{D},i} = 0
	\end{cases}
\end{equation}
Furthermore, we add the breakpoint $t_{-1}$ corresponding to maximum allowed stepsize on $\mathcal{M}$, which is usually taken as the injectivity radius of the exponential map of $\mathcal{M}$ at $p_{\mathrm{M}}$.
Next, we sort the breakpoints in nondecreasing order and denote the sorted sequence by $t_{(1)} \leq t_{(2)} \leq \cdots \leq t_{(n+1)}$.

In the next step we calculate the initial linear and quadratic terms of the piecewise quadratic function $q(t) = m(d_{\mathrm{PL}}(t)) = q(0) + f' t + \frac{1}{2}f'' t^2$ on the first segment $[0, t_{(1)}]$, $f' = \langle \operatorname{grad} f(p), d \rangle_{p}$ and $f'' = \langle d, H[d] \rangle_{p}$.
If either $f' = 0$ or $f'' = 0$, then we return the zero tangent vector and the status \texttt{NOT FOUND}.
Otherwise, we compute the minimizer of $q(t)$ on the first segment, $\Delta t_{\mathrm{min}} = -\frac{f'}{f''}$.

Next, suppose we are examining the interval $[t_{(j)}, t_{(j+1)}]$ for some $j \in \{1, 2, \dots, n\}$.
Denote $Z_{(j)} = d_{\mathrm{PL}}(t_{(j)})$,  $\Delta t = t - t_{(j)}$, $\Delta t_{(j)} = t_{(j+1)} - t_{(j)}$ and $\hat{d}_{(j)} = (\hat{d}_{\mathrm{D}}, d_{\mathcal{M}})$ where
\begin{equation}
	\hat{d}_{\mathrm{D},i} = \begin{cases}
		d_{\mathrm{D},i} & \text{ when } t_{(j)} < t_i \\
		0 & \text{ otherwise, }
	\end{cases}
\end{equation}
for each $i = 1, 2, \dots, n$.
Now we can write
\begin{align*}
	q(t_{(j)} + \Delta t) = &  m(Z_{(j)} + \Delta t \hat{d}_{(j)}) \\
	= & f(p) + \langle g, Z_{(j)} + \Delta t \hat{d}_{(j)} \rangle_{p} + \frac{1}{2} \langle Z_{(j)} + \Delta t \hat{d}_{(j)}, H[Z_{(j)} + \Delta t \hat{d}_{(j)}] \rangle_{p}\\
	= & f(p) + \langle g, Z_{(j)} \rangle_{p} + \frac{1}{2} \langle Z_{(j)}, H[Z_{(j)}] \rangle_{p} \\
	& + \Delta t \left( \langle g, \hat{d}_{(j)} \rangle_{p} + \langle \hat{d}_{(j)}, H[Z_{(j)}] \rangle_{p} \right) \\
	& + \frac{1}{2} \Delta t^2 \langle \hat{d}_{(j)}, H[\hat{d}_{(j)}] \rangle_{p}.
\end{align*}
The quadratic surrogate $\hat{m}(\Delta t)$ on the current segment is then given by
\begin{equation}
	\hat{m}(\Delta t) = q(t_{(j)} + \Delta t) = f_{(j)} + f_{(j)}' \Delta t + \frac{1}{2} f_{(j)}'' \Delta t^2,
\end{equation}
where
\begin{align*}
	f_{(j)} & = f(p) + \langle g, Z_{(j)} \rangle_{p} + \frac{1}{2} \langle Z_{(j)}, H[Z_{(j)}] \rangle_{p}, \\
	f_{(j)}' & = \langle g, \hat{d}_{(j)} \rangle_{p} + \langle \hat{d}_{(j)}, H[Z_{(j)}] \rangle_{p}, \\
	f_{(j)}'' & = \langle \hat{d}_{(j)}, H[\hat{d}_{(j)}] \rangle_{p}.
\end{align*}
Since we assume that $H$ is positive definite $\hat{m}$ has a minimizer $\Delta t_{\mathrm{min}} = -f_{(j)}' / f_{(j)}''$.
If $t_{(j)} + \Delta t_{\mathrm{min}}$ lies within the current segment, then we have found the first minimizer of $q(t)$.
If $t_{(j)} + \Delta t_{\mathrm{min}} < t_{(j)}$, then the minimum lies at $t_{(j)}$.
Otherwise, if $t_{(j)} + \Delta t_{\mathrm{min}} > t_{(j+1)}$, then we need to proceed to the next segment.

When proceeding to the next segment, we don't have to recompute $f_{(j+1)}'$ and $f_{(j+1)}''$ from scratch.
We can update them from $f_{(j)}'$ and $f_{(j)}''$.
Let $b$ be the index of the component that reached its bound at $t_{(j+1)}$.
Then, we have $\hat{d}_{(j+1)} = \hat{d}_{(j)} - d_{\mathrm{D},b} e_b$ and $Z_{(j+1)} = Z_{(j)} + \Delta t_{(j)} \hat{d}_{(j)}$.
Using these relations, we can derive the following update formulas:
\begin{align*}
	f_{(j+1)}' =& \langle g, \hat{d}_{(j+1)} \rangle_{p} + \langle \hat{d}_{(j+1)}, H[Z_{(j+1)}] \rangle_{p}\\
	=& \langle g, \hat{d}_{(j)} - d_{\mathrm{D},b} e_b \rangle_{p} + \langle \hat{d}_{(j)} - d_{\mathrm{D},b} e_b, H[Z_{(j)} + \Delta t_{(j)} \hat{d}_{(j)}] \rangle_{p}\\
	=& \langle g, \hat{d}_{(j)} \rangle_{p} - d_{\mathrm{D},b} \langle g, e_b \rangle_{p}  + \langle \hat{d}_{(j)}, H[Z_{(j)}] \rangle_{p} + \Delta t_{(j)} \langle \hat{d}_{(j)}, H[\hat{d}_{(j)}] \rangle_{p} - \\
	& \quad d_{\mathrm{D},b} \langle e_b, H[Z_{(j)}] \rangle_{p} - d_{\mathrm{D},b} \Delta t_{(j)} \langle e_b, H[\hat{d}_{(j)}] \rangle_{p} \\
	= & f_{(j)}' - d_{\mathrm{D},b} g_{\mathrm{D},b} + \Delta t_{(j)} f_{(j)}'' - d_{\mathrm{D},b} \langle e_b, H[Z_{(j+1)}] \rangle_{p}.
\end{align*}
Similarly, we derive
\begin{align*}
	f_{(j+1)}'' =& \langle \hat{d}_{(j+1)}, H[\hat{d}_{(j+1)}] \rangle_{p} \\
	=& \langle \hat{d}_{(j)} - d_{\mathrm{D},b} e_b, H[\hat{d}_{(j)} - d_{\mathrm{D},b} e_b]\rangle_{p} \\
	=& \langle \hat{d}_{(j)}, H[\hat{d}_{(j)}] \rangle_{p} - 2 d_{\mathrm{D},b} \langle e_b, H[\hat{d}_{(j)}] \rangle_{p} + d_{\mathrm{D},b}^2 \langle e_b, H[e_b] \rangle_{p} \\
	= & f_{(j)}'' - 2 d_{\mathrm{D},b} \langle e_b, H[\hat{d}_{(j)}] \rangle_{p} + d_{\mathrm{D},b}^2 \langle e_b, H[e_b] \rangle_{p}.
\end{align*}

Now, there are exactly three Hessian values that need to be computed in each segment transition: $\langle e_b, H[e_b] \rangle_{p}$, $\langle e_b, H[\hat{d}_{(j)}] \rangle_{p}$ and $\langle e_b, H[Z_{(j+1)}] \rangle_{p}$.
The first one can be computed statelessly in the selected Hessian approximation.
The other two values require maintaining certain state between segment transitions to avoid recomputing them from scratch.

Next, we prepare for iteration over the segments defined by the breakpoints.
We denote:
\begin{align*}
	& p_{\mathrm{y},(j)} = [\langle y_{(i)}, \hat{d}_{(j)} \rangle_{p}]_{i=1, 2, \dots, \mu}, \\
	& p_{\mathrm{s},(j)} = [\theta \langle s_{(i)}, \hat{d}_{(j)} \rangle_{p}]_{i=1, 2, \dots, \mu}, \\
	& c_{\mathrm{y},(j)} = [\langle y_{(i)}, Z_{(j)} \rangle_{p}]_{i=1, 2, \dots, \mu}, \\
	& c_{\mathrm{s},(j)} = [\theta \langle s_{(i)}, Z_{(j)} \rangle_{p}]_{i=1, 2, \dots, \mu}.
\end{align*}
Given these vectors, we can compute the required Hessian values as
\begin{align*}
	\langle e_b, H[\hat{d}_{(j)}] \rangle_{p} & = H_{\mathrm{W}}(\theta d_b, e_{\mathrm{D},b}, e_{\mathrm{D},b},  p_{\mathrm{y},(j)}, p_{\mathrm{s},(j)}), \\
	\langle e_b, H[Z_{(j)}] \rangle_{p} & = H_{\mathrm{W}}(\theta t_{(j)}, e_{\mathrm{D},b}, e_{\mathrm{D},b},  c_{\mathrm{y},(j)}, c_{\mathrm{s},(j)}),
\end{align*}
see Eq.~\eqref{eq:H_W_k} for the definition of $H_{\mathrm{W}}$.
Now, we can initialize the vectors at the start of the first segment as
\begin{align*}
	p_{\mathrm{y},(1)} & = [\langle y_{(i)}, d \rangle_{p}]_{i=1, 2, \dots, \mu}, \\
	p_{\mathrm{s},(1)} & = [\theta \langle s_{(i)}, d \rangle_{p}]_{i=1, 2, \dots, \mu}, \\
	c_{\mathrm{y},(1)} & = [0, 0, \dots, 0]^{\top} \in \mathbb{R}^{\mu}, \\
	c_{\mathrm{s},(1)} & = [0, 0, \dots, 0]^{\top} \in \mathbb{R}^{\mu}.
\end{align*}
Finally, we can update the vectors $c_{\mathrm{y},(j)}$, $c_{\mathrm{s},(j)}$, $p_{\mathrm{y},(j)}$ and $p_{\mathrm{s},(j)}$ between segments as follows:
\begin{align*}
	p_{\mathrm{y},(j+1)} & = p_{\mathrm{y},(j)} - d_b [y_{\mathrm{D},i,b}]_{i=1, 2, \dots, \mu}, \\
	p_{\mathrm{s},(j+1)} & = p_{\mathrm{s},(j)} - d_b \theta [s_{\mathrm{D},i,b}]_{i=1, 2, \dots, \mu}, \\
	c_{\mathrm{y},(j+1)} & = c_{\mathrm{y},(j)} + \Delta t_{(j)} p_{\mathrm{y},(j)}, \\
	c_{\mathrm{s},(j+1)} & = c_{\mathrm{s},(j)} + \Delta t_{(j)} p_{\mathrm{s},(j)}.
\end{align*}

The part of the algorithm that depends on the Hessian approximation is encapsulated within functions $\texttt{quadratic\_segment\_surrogate}$ (Algorithm~\ref{alg:qss_init}) and $\texttt{hessian\_segment\_values}$ (Algorithm~\ref{alg:qss_hessian_values}).
It is defined by an implicit vector $d_z$ tangent to $p$, pointing at the end of the current segment.
The function $\texttt{quadratic\_segment\_surrogate}$ initializes a temporary storage structure $q_s$ that is used to compute Hessian values in subsequent segments.
The storage corresponds to variables $p$ and $c$ in the L-BFGS-B implementation from~\cite{ByrdLuNocedalZhu:1995}.
The function $\texttt{hessian\_segment\_values}$ computes the required Hessian values in the current segment using $q_s$.

\begin{algorithm}
\caption{Quadratic segment surrogate initialization (L-BFGS-B variant)}
\label{alg:qss_init}
\begin{algorithmic}[1]
	\Require Current iterate $p$, search direction $d = -B \operatorname{grad} f(p)$, Hessian approximation $H$
	\Ensure Temporary storage structure $q_s$ required for calculation of Hessian values in subsequent segments
	\State $c_s = [0, 0, \dots, 0]^{\top} \in \mathbb{R}^{\mu}$
	\State $c_y = [0, 0, \dots, 0]^{\top} \in \mathbb{R}^{\mu}$
	\State $p_s = \theta [\langle s_i, d \rangle_p]_{i=1, 2, \dots, \mu}$
	\State $p_y = [\langle y_i, d \rangle_p]_{i=1, 2, \dots, \mu}$
	\State \Return $q_s = (c_s, c_y, p_s, p_y)$
\end{algorithmic}
\end{algorithm}

\begin{algorithm}
\caption{Hessian calculation for subsequent segments (L-BFGS-B variant). $M$ is the precomputed matrix from Eq.~\eqref{eq:M_k} for the current iteration.}
\label{alg:qss_hessian_values}
\begin{algorithmic}[1]
	\Require Temporary storage $q_s$, current iterate $p$, current step size $t$, length of the current segment $\Delta t$, current index $b$, value of direction vector at index $b$, Hessian approximation $H$
	\Ensure $v_1 = \langle e_b, H[Z] \rangle$, $v_2 = \langle e_b, H[\hat{d}] \rangle$

	\State $c_s, c_y, p_s, p_y \gets q_s$
	\State $c_y \gets c_y + \Delta t \cdot p_y$
	\State $c_s \gets c_s + \Delta t \cdot p_s$
	\State $\xi_y \gets [y_{\mathrm{D},i,b} ]_{i=1, 2, \dots, \mu}$
	\State $\xi_s \gets \theta [s_{\mathrm{D},i,b} ]_{i=1, 2, \dots, \mu}$
	\State $v_1 \gets \theta t d_b + \begin{bmatrix}
		\xi_y^{\top} & \xi_s^{\top}
	\end{bmatrix} M \begin{bmatrix}
		c_y \\
		c_s
	\end{bmatrix}$
	\State $v_2 \gets \theta d_b + \begin{bmatrix}
		\xi_y^{\top} & \xi_s^{\top}
	\end{bmatrix} M \begin{bmatrix}
		p_y \\
		p_s
	\end{bmatrix}$
	\State $p_y \gets p_y - d_b \xi_y$
	\State $p_s \gets p_s - d_b \xi_s$
	\State $q_s \gets (c_s, c_y, p_s, p_y)$
	\State \Return $(v_1, v_2)$
\end{algorithmic}
\end{algorithm}

Evaluation of subsequent segments proceeds similarly to the original generalized Cauchy point algorithm from~\cite{ByrdLuNocedalZhu:1995}.
The loop has an additional stopping criterion that checks whether we have reached the maximum allowed stepsize on $\mathcal{M}$.
After the segment with the minimizer is found, we construct the direction $d_{\mathrm{GCD}}$ by multiplying the original direction by the time step at which the minimum is found and fixing the components that reached their bounds at the corresponding breakpoints.
Additionally, the maximum stepsize for the subsequent line search is also computed.
The entire procedure is described by Algorithm~\ref{alg:generalized_cauchy_direction}.

\begin{algorithm}
\caption{Generalized Cauchy direction computation}
\label{alg:generalized_cauchy_direction}
\begin{algorithmic}[1]
	\Require Current iterate $p$, gradient $X = \operatorname{grad} f(p)$, search direction $d = -B \operatorname{grad} f(p)$, Hessian approximation $H$, maximum stepsize $t_{\mathcal{M},\mathrm{max}}$ on $\mathcal{M}$
	\Ensure New search direction $d_{\mathrm{GCD}}$
	\Ensure Search result $r$ (one of : $\texttt{FOUND LIMITED}, \texttt{FOUND UNLIMITED}, \texttt{NOT FOUND}$)
	\Ensure Maximum stepsize $t_{\mathrm{GCD},\mathrm{max}}$ for the subsequent line search, or $-1$ if the result is \texttt{NOT FOUND}
	\State $\mathcal{F} \gets \{(t_{(i)}, i) : t_{(i)} > 0\}$ \Comment{min-heap of breakpoints}
	\If{there is a finite breakpoint in $\mathcal{F}$}
		\State $f_{\mathrm{FL}} \gets $ \textbf{true}
	\Else
		\State $f_{\mathrm{FL}} \gets $ \textbf{false}
	\EndIf
	\State \Call{push}{$\mathcal{F}$, $(t_{\mathcal{M},\mathrm{max}}, -1)$}
	\State $f' \gets \langle \operatorname{grad} f(p), d \rangle_{p}$ \Comment{linear term of $q(t)$ in the current segment}
	\State $f'' \gets \langle d, H[d] \rangle_{p}$ \Comment{quadratic term of $q(t)$ in the current segment}
	\If{$f' = 0$ or $f'' = 0$}
		\State \Return $0 \cdot d$, \texttt{NOT FOUND}, $-1$
	\EndIf
	\State $\Delta t_{\mathrm{min}} \gets -\frac{f'}{f''}$
	\State $t_{\mathrm{old}} \gets 0$
	\State $t, b \gets$ \Call{pop\_min}{$\mathcal{F}$} \Comment{distance to and index of the next breakpoint}
	\State $\Delta t \gets t$
	\State $q_s \gets$ \Call{quadratic\_segment\_surrogate}{$p$, $d$, $H$}
	\While{$\Delta t_{\mathrm{min}} > \Delta t$ and $b \neq -1$}  \Comment{examine the next segment}
		\State $d_b \gets d_{\mathrm{D},b}$
		\State $g_b \gets X_{\mathrm{D},b}$
		\State $v_1, v_2 \gets $ \Call{hessian\_segment\_values}{$q_s$, $p$, $t$, $\Delta t$, $b$, $d_b$, $H$} \Comment{Hessian values in the current segment}
		\State $f' \gets f' + \Delta t \cdot f'' - d_b (g_b + v_1)$ \Comment{update linear term of $q(t)$}
		\State $f'' \gets f'' - 2 d_b v_2 + d_b^2 \langle e_b, H[e_b] \rangle$ \Comment{update quadratic term of $q(t)$}
		\State $t_{\mathrm{old}} \gets t$
		\If{$f'=0$ or $f''=0$}
			\State $\Delta t_{\mathrm{min}} \gets 0$
			\State \textbf{break}
		\EndIf
		\State $\Delta t_{\mathrm{min}} \gets -\frac{f'}{f''}$
		\If{\Call{is\_empty}{$\mathcal{F}$}}
			\State \textbf{break}
		\EndIf
		\State $t, b \gets$ \Call{pop\_min}{$\mathcal{F}$}
		\State $\Delta t \gets t - t_{\mathrm{old}}$
	\EndWhile
	\State $\Delta t_{\mathrm{min}} \gets \max(0, \Delta t_{\mathrm{min}})$
	\State $t_{\mathrm{old}} \gets t_{\mathrm{old}} + \Delta t_{\mathrm{min}}$ \Comment{Beginning of the last searched segment plus distance to the minimizer in this segment}
	\State $d_{\mathrm{GCD}} \gets t_{\mathrm{old}} \cdot d$ \Comment{initialize GCD direction}
	\For{$i = 1, 2, \dots, n$}
		\If{$t_{(i)} < t$}\Comment{the $i$th component was fixed at a bound}
			\If{$d_{\mathrm{GCD},\mathrm{D},i} > 0$}
				\State $d_{\mathrm{GCD},\mathrm{D},i} \gets u_i - p_{\mathrm{D},i}$
			\Else
				\State $d_{\mathrm{GCD},\mathrm{D},i} \gets l_i - p_{\mathrm{D},i}$
			\EndIf
		\EndIf
	\EndFor
	\If{$f_{\mathrm{FL}}$}
		\State $t_{\mathrm{ni}} \gets \min \{t_{\mathcal{M}, \mathrm{max}}\} \cup \{ t_{(i)} : t_{(i)} > 0 \}$ \Comment{distance to the nearest inactive bound or maximum stepsize on $\mathcal{M}$}
		\State \Return $d_{\mathrm{GCD}}$, \texttt{FOUND LIMITED}, $\max(1, t_{\mathrm{ni}} / t_{\mathrm{old}})$
	\Else
		\State \Return $d_{\mathrm{GCD}}$, \texttt{FOUND UNLIMITED}, $\infty$
	\EndIf

\end{algorithmic}
\end{algorithm}

\section{Results}
\label{sec:results}

\subsection{Euclidean test problems with bound constraints}

Comparison with existing Euclidean solvers on problems with bounds constraints was performed to establish how much performance is lost by extending the domain from a pure hypercube $D$ to the product $D \times \mathcal{M}$.
Problems from the CUTEst set~\cite{GouldOrbanToint:2015,GrattonToint:2024} with only bounds constraints were selected for this test.
IPOPT~\cite{WachterBiegler:2006}, L-BFGS-B~\cite{ByrdLuNocedalZhu:1995} and UNO solver~\cite{VanaretLeyffer:2026} were initially considered as Euclidean baselines.
Since the standard L-BFGS-B achieved the best performance on the selected problems, it was selected as a baseline for further analysis.
UNO solver 0.3.0, IPOPT.jl 1.14.1 and LBFGSB.jl 0.4.1 were used for the experiments.
The forthcoming version 0.6.0 of Manopt.jl was used for Riemannian L-BFGS-B implementation.
Experiments were performed on a machine running Linux Mint 22.3 on an AMD Ryzen 9 9950X3D CPU.

\begin{figure}
\begin{centering}

\begin{tikzpicture}
\begin{axis}[
  width=12cm,
  height=8cm,
  xlabel={Performance ratio},
  ylabel={Fraction of problems solved},
  xmin=0,
  xmax=13.019381125380816,
  ymin=0,
  ymax=1,
  grid=major,
  legend pos=south east,
  legend cell align={left},
  xticklabel={$2^{\pgfmathprintnumber{\tick}}$},
]
\addplot+[const plot, no markers, thick, color=blue, solid] coordinates {(0.0,0.43333333333333335) (0.0,0.43333333333333335) (0.05094945517593566,0.4666666666666667) (0.15541314381932578,0.4666666666666667) (0.16150834450860105,0.4666666666666667) (0.17263978916481262,0.5) (0.17572919993722164,0.5) (0.18492391142850773,0.5) (0.2007923742907667,0.5333333333333333) (0.2116401603124625,0.5666666666666667) (0.2942577552063276,0.6) (0.33666108940440337,0.6) (0.34473269666148804,0.6333333333333333) (0.35041084684693335,0.6333333333333333) (0.36887321670172635,0.6333333333333333) (0.37773007829176003,0.6333333333333333) (0.40471710721860893,0.6666666666666666) (0.4242496900824379,0.7) (0.45596532061962225,0.7) (0.4766517818304484,0.7333333333333333) (0.5752836158253216,0.7666666666666667) (0.6607527410923582,0.7666666666666667) (1.0350022092752011,0.7666666666666667) (1.351817011973616,0.8) (1.3783777489393265,0.8333333333333334) (1.4015623782672841,0.8333333333333334) (1.5504614034474116,0.8666666666666667) (1.675826644423905,0.8666666666666667) (2.4869037956758158,0.9) (3.041243451784948,0.9) (3.7508201611662044,0.9333333333333333) (4.678946661586113,0.9333333333333333) (4.749246361435427,0.9333333333333333) (4.89705977599518,0.9333333333333333) (5.227468559545668,0.9333333333333333) (5.268949001485523,0.9333333333333333) (5.314746494753967,0.9333333333333333) (5.653381192528148,0.9333333333333333) (5.796751368139142,0.9333333333333333) (5.849133115019879,0.9333333333333333) (5.8862044550979995,0.9333333333333333) (5.944954173679409,0.9333333333333333) (6.116174581870618,0.9333333333333333) (6.194665339669838,0.9333333333333333) (6.218033166525744,0.9333333333333333) (7.003230089054143,0.9333333333333333) (7.124320940401024,0.9666666666666667) (7.251908047747934,0.9666666666666667) (7.255826177867817,0.9666666666666667) (7.593292251412594,0.9666666666666667) (7.827227323673752,0.9666666666666667) (8.198983957637084,0.9666666666666667) (8.557327834182082,0.9666666666666667) (8.71282197510104,0.9666666666666667) (8.843258862159896,1.0) (11.134473938614128,1.0) (11.20570098459994,1.0) (11.253949134491137,1.0) (12.254512289912114,1.0) (12.669227492221832,1.0) (12.754201418978411,1.0) (13.019381125380816,1.0)};
\addlegendentry{Manopt.jl}
\addplot+[const plot, no markers, thick, color=red, solid] coordinates {(0.0,0.5333333333333333) (0.0,0.5333333333333333) (0.05094945517593566,0.5333333333333333) (0.15541314381932578,0.5666666666666667) (0.16150834450860105,0.6) (0.17263978916481262,0.6) (0.17572919993722164,0.6333333333333333) (0.18492391142850773,0.6666666666666666) (0.2007923742907667,0.6666666666666666) (0.2116401603124625,0.6666666666666666) (0.2942577552063276,0.6666666666666666) (0.33666108940440337,0.7) (0.34473269666148804,0.7) (0.35041084684693335,0.7333333333333333) (0.36887321670172635,0.7666666666666667) (0.37773007829176003,0.8) (0.40471710721860893,0.8) (0.4242496900824379,0.8) (0.45596532061962225,0.8333333333333334) (0.4766517818304484,0.8333333333333334) (0.5752836158253216,0.8333333333333334) (0.6607527410923582,0.8666666666666667) (1.0350022092752011,0.9) (1.351817011973616,0.9) (1.3783777489393265,0.9) (1.4015623782672841,0.9) (1.5504614034474116,0.9) (1.675826644423905,0.9333333333333333) (2.4869037956758158,0.9333333333333333) (3.041243451784948,0.9666666666666667) (3.7508201611662044,0.9666666666666667) (4.678946661586113,0.9666666666666667) (4.749246361435427,0.9666666666666667) (4.89705977599518,0.9666666666666667) (5.227468559545668,1.0) (5.268949001485523,1.0) (5.314746494753967,1.0) (5.653381192528148,1.0) (5.796751368139142,1.0) (5.849133115019879,1.0) (5.8862044550979995,1.0) (5.944954173679409,1.0) (6.116174581870618,1.0) (6.194665339669838,1.0) (6.218033166525744,1.0) (7.003230089054143,1.0) (7.124320940401024,1.0) (7.251908047747934,1.0) (7.255826177867817,1.0) (7.593292251412594,1.0) (7.827227323673752,1.0) (8.198983957637084,1.0) (8.557327834182082,1.0) (8.71282197510104,1.0) (8.843258862159896,1.0) (11.134473938614128,1.0) (11.20570098459994,1.0) (11.253949134491137,1.0) (12.254512289912114,1.0) (12.669227492221832,1.0) (12.754201418978411,1.0) (13.019381125380816,1.0)};
\addlegendentry{L-BFGS-B}
\addplot+[const plot, no markers, thick, color=teal, solid] coordinates {(0.0,0.03333333333333333) (0.0,0.03333333333333333) (0.05094945517593566,0.03333333333333333) (0.15541314381932578,0.03333333333333333) (0.16150834450860105,0.03333333333333333) (0.17263978916481262,0.03333333333333333) (0.17572919993722164,0.03333333333333333) (0.18492391142850773,0.03333333333333333) (0.2007923742907667,0.03333333333333333) (0.2116401603124625,0.03333333333333333) (0.2942577552063276,0.03333333333333333) (0.33666108940440337,0.03333333333333333) (0.34473269666148804,0.03333333333333333) (0.35041084684693335,0.03333333333333333) (0.36887321670172635,0.03333333333333333) (0.37773007829176003,0.03333333333333333) (0.40471710721860893,0.03333333333333333) (0.4242496900824379,0.03333333333333333) (0.45596532061962225,0.03333333333333333) (0.4766517818304484,0.03333333333333333) (0.5752836158253216,0.03333333333333333) (0.6607527410923582,0.03333333333333333) (1.0350022092752011,0.03333333333333333) (1.351817011973616,0.03333333333333333) (1.3783777489393265,0.03333333333333333) (1.4015623782672841,0.06666666666666667) (1.5504614034474116,0.06666666666666667) (1.675826644423905,0.06666666666666667) (2.4869037956758158,0.06666666666666667) (3.041243451784948,0.06666666666666667) (3.7508201611662044,0.06666666666666667) (4.678946661586113,0.1) (4.749246361435427,0.13333333333333333) (4.89705977599518,0.16666666666666666) (5.227468559545668,0.16666666666666666) (5.268949001485523,0.2) (5.314746494753967,0.23333333333333334) (5.653381192528148,0.26666666666666666) (5.796751368139142,0.3) (5.849133115019879,0.3333333333333333) (5.8862044550979995,0.36666666666666664) (5.944954173679409,0.4) (6.116174581870618,0.43333333333333335) (6.194665339669838,0.4666666666666667) (6.218033166525744,0.5) (7.003230089054143,0.5333333333333333) (7.124320940401024,0.5333333333333333) (7.251908047747934,0.5666666666666667) (7.255826177867817,0.6) (7.593292251412594,0.6333333333333333) (7.827227323673752,0.6666666666666666) (8.198983957637084,0.7) (8.557327834182082,0.7333333333333333) (8.71282197510104,0.7666666666666667) (8.843258862159896,0.7666666666666667) (11.134473938614128,0.8) (11.20570098459994,0.8333333333333334) (11.253949134491137,0.8666666666666667) (12.254512289912114,0.9) (12.669227492221832,0.9333333333333333) (12.754201418978411,0.9666666666666667) (13.019381125380816,1.0)};
\addlegendentry{UNO Solver}
\end{axis}
\end{tikzpicture}

\end{centering}
\label{fig:performance_profile_euclidean}
\caption{Performance plot}
\end{figure}
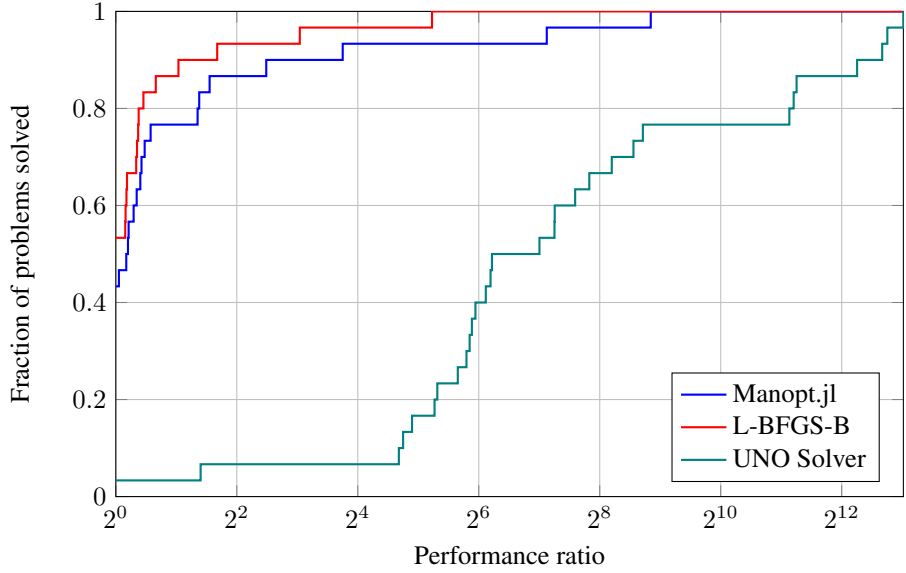

In the experiments the standard termination conditions of L-BFGS-B with the relative a posteriori cost change factor parameter equal to 1000, maximum number of iterations equal to 1000 and projected gradient norm tolerance parameter equal to $10^{-12}$.
Summarized results are shown in Figure~\ref{fig:performance_profile_euclidean}, while details are shown in Table~\ref{tab:euclidean_results_all}.
Riemannian L-BFGS-B is nearly on par with the standard L-BFGS-B.
The UNO solver is significantly slower on most problems than both L-BFGS-B and Riemannian L-BFGS-B.

\begin{sidewaystable}
\centering
	\caption{Detailed results for Euclidean problems with bound constraints
	comparing our solver from Manopt.jl (M) with the Fortran L-BFGS-B.jl (F)
	and the UNO solver (U).
	Best values for each problem are shown in bold.}
	\label{tab:euclidean_results_all}
\begin{tabular}{l
  S[table-format=4.5]
  S[table-format=2.5]
  S[table-format=3.3]
  S[table-format=3.0]
  S[table-format=3.0]
  S[table-format=4.0]
  S[table-format=3.0]
  S[table-format=3.0]
  S[table-format=4.0]
  S[table-format=1.3e-2,tight-spacing=true,scientific-notation = true,]
  S[table-format=1.3e-2,tight-spacing=true,scientific-notation = true,]
  S[table-format=1.3e-2,tight-spacing=true,scientific-notation = true,]}
\toprule
{\bfseries problem} & \multicolumn{3}{c}{{\bfseries time} [ms]} & \multicolumn{3}{c}{\bfseries objective calls} & \multicolumn{3}{c}{\bfseries gradient calls} & \multicolumn{3}{c}{\bfseries objective value} \\
 & M & F & U & M & F & U & M & F & U & M & F & U \\
\cmidrule(rl){1-1}\cmidrule(r){2-4}\cmidrule(lr){5-7}\cmidrule(lr){8-10}\cmidrule(l){11-13}
ALLINIT & \bfseries 0.04027 & 0.05524 & 2.481 & 28 & \bfseries 18 & 51 & \bfseries 18 & \bfseries 18 & 51 & 16.71 & \bfseries 16.71 & \bfseries 16.71 \\
BDEXP & \bfseries 11.72 & 15.13 & 651.3 & 56 & \bfseries 41 & 45 & \bfseries 30 & 41 & 45 & 2.606e-10 & \bfseries 1.336e-10 & 8.471e-07 \\
BRANIN & 0.03293 & \bfseries 0.02844 & 1.132 & 24 & \bfseries 14 & 20 & 16 & \bfseries 14 & 20 & \bfseries 0.3979 & \bfseries 0.3979 & \bfseries 0.3979 \\
CAMEL6 & \bfseries 0.0346 & 0.03909 & 0.8864 & 21 & \bfseries 16 & 17 & \bfseries 13 & 16 & 17 & \bfseries -1.032 & \bfseries -1.032 & \bfseries -1.032 \\
CHARDIS0 & 113.3 & \bfseries 81.39 & 1.922e+05 & 10 & \bfseries 5 & 13908 & 7 & \bfseries 5 & 13908 & \bfseries 0 & \bfseries 1.227e-28 & 2.399e-10 \\
\midrule
CHARDIS02 & 6918 & \bfseries 49.58 & 1.115e+05 & 678 & \bfseries 3 & 6445 & 369 & \bfseries 3 & 6445 & \bfseries 1.993e+06 & 1.466e+09 & 3.492e+06 \\
DGOSPEC & 0.05537 & \bfseries 0.0189 & 5.555 & 33 & \bfseries 12 & 116 & 18 & \bfseries 12 & 116 & -997 & -997 & \bfseries -997 \\
EG1 & 0.04973 & \bfseries 0.03756 & 1.891 & 35 & \bfseries 12 & 38 & 21 & \bfseries 12 & 38 & -1.132 & -1.133 & \bfseries -1.133 \\
EXPLIN & \bfseries 0.5643 & 4.646 & 3676 & \bfseries 22 & 177 & 5500 & \bfseries 21 & 177 & 5500 & -6.986e+07 & \bfseries -7.193e+07 & -7.192e+07 \\
EXPLIN2 & \bfseries 0.5551 & 1.774 & 3835 & \bfseries 22 & 66 & 5583 & \bfseries 21 & 66 & 5583 & -7.183e+07 & -7.2e+07 & \bfseries -7.2e+07 \\
\midrule
EXPQUAD & \bfseries 9.054 & 11.76 & 3799 & 309 & \bfseries 213 & 5461 & \bfseries 177 & 213 & 5460 & -3.685e+09 & \bfseries -3.685e+09 & -3.649e+09 \\
GRIDGENA & 71.83 & \bfseries 5.336 & 2.607e+04 & 248 & \bfseries 21 & 11732 & 135 & \bfseries 21 & 11732 & \bfseries 2.352e+04 & {NaN} & 4.499e+04 \\
HADAMALS & \bfseries 0.9433 & 1.051 & 355.4 & 23 & \bfseries 20 & 654 & \bfseries 14 & 20 & 654 & 7312 & 7312 & \bfseries 6412 \\
HART6 & 0.06753 & \bfseries 0.05876 & 4.304 & 37 & \bfseries 21 & 79 & \bfseries 21 & \bfseries 21 & 79 & -3.323 & -3.323 & \bfseries -3.323 \\
HIMMELP1 & \bfseries 0.01785 & 0.02254 & 1.329 & 18 & \bfseries 13 & 32 & 14 & \bfseries 13 & 32 & -23.9 & -23.9 & \bfseries -62.05 \\
\midrule
HOLMES & 75.82 & \bfseries 50.89 & 134.4 & 189 & \bfseries 84 & 144 & 106 & \bfseries 84 & 144 & 1249 & 1248 & \bfseries 1248 \\
HS38 & \bfseries 0.09975 & 0.1116 & 15.2 & 33 & \bfseries 27 & 489 & \bfseries 18 & 27 & 489 & \bfseries 2.628e-18 & \bfseries 5.221e-17 & \bfseries 5.765e-18 \\
HS4 & 0.00378 & \bfseries 0.003649 & 0.5577 & \bfseries 2 & \bfseries 2 & 9 & \bfseries 2 & \bfseries 2 & 9 & 2.667 & 2.667 & \bfseries 2.667 \\
HS45 & \bfseries 0.008226 & 0.01049 & 1.868 & \bfseries 10 & 11 & 45 & \bfseries 8 & 11 & 45 & 1 & 1 & \bfseries 1 \\
HS5 & 0.02723 & \bfseries 0.02029 & 0.7825 & 18 & \bfseries 9 & 12 & 12 & \bfseries 9 & 12 & \bfseries -1.913 & \bfseries -1.913 & \bfseries -1.913 \\
\midrule
LOGROS & 0.35 & \bfseries 0.2756 & 16.3 & 296 & \bfseries 112 & 1067 & 208 & \bfseries 112 & 1067 & \bfseries 0 & \bfseries 8.882e-16 & \bfseries 0 \\
MAXLIKA & \bfseries 0.7117 & 26.67 & 91.31 & \bfseries 44 & 1038 & 1613 & \bfseries 24 & 1038 & 1613 & 1150 & 1137 & \bfseries 1136 \\
MCCORMCK & 14.41 & \bfseries 5.647 & 325.5 & 61 & \bfseries 14 & 28 & 44 & \bfseries 14 & 28 & -4567 & -4567 & \bfseries -4567 \\
MDHOLE & 0.2427 & \bfseries 0.2153 & 6.416 & 181 & \bfseries 89 & 245 & 124 & \bfseries 89 & 245 & 9.436e-28 & 1.43e-33 & \bfseries -1e-08 \\
MINSURFO & \bfseries 158.5 & 324.8 & 1.099e+04 & \bfseries 437 & 544 & 709 & \bfseries 220 & 544 & 709 & 2.507 & \bfseries 2.507 & 2.507 \\
\midrule
POWELLBC & 1809 & 6.225 & \bfseries 3.938 & 1288 & 3 & \bfseries 2 & 696 & 3 & \bfseries 0 & \bfseries 3.106e+05 & 2.323e+06 & {Inf} \\
PROBPENL & 0.1115 & \bfseries 0.0909 & 754.8 & 5 & \bfseries 4 & 638 & \bfseries 4 & \bfseries 4 & 638 & 3.992e-07 & 3.992e-07 & \bfseries -2.228e-05 \\
QRTQUAD & \bfseries 69 & 78.44 & 1.332e+04 & 583 & \bfseries 524 & 2159 & \bfseries 338 & 524 & 2159 & -1.218e+09 & \bfseries -1.218e+09 & -1.217e+09 \\
S368 & 0.3501 & \bfseries 0.06245 & 1.68 & 143 & \bfseries 13 & 22 & 117 & \bfseries 13 & 22 & -0.75 & -0.75 & \bfseries -0.75 \\
SINEALI & 3.733 & \bfseries 1.436 & 3507 & 104 & \bfseries 24 & 8754 & 71 & \bfseries 24 & 8754 & \bfseries -9.99e+04 & -9.987e+04 & -9.972e+04 \\
\bottomrule
\end{tabular}
\end{sidewaystable}

\subsection{Amplitude-limited blind source separation with Gaussianity penalization}

In the problem of blind source separation we have $k$ independent signal sources of length $n$ that are linearly mixed together, and we want to recover them from the data observed using $r$ sensors.
The observed data is collected in a matrix $X \in \mathbb{R}^{r \times n}$.
Here we assume that the sources are bounded in amplitude, i.e., $S \in [-A, A]^n$ for some $A > 0$.
Additionally, we want to promote the non-Gaussianity of the recovered sources, which is a common approach in blind source separation problems.
Here, we take an approach based on negentropy maximization through~\cite{HyvarinenOja:2000}, which leads to the following optimization problem:
\begin{equation}
	\operatorname*{arg\ min}_{W\in \mathrm{St}(k, r) , S \in [-A, A]^n}  \frac{1}{2}\lVert S - W X \rVert^2 - \lambda \sum_{i=1}^n \log(\cosh(S_i)),
\end{equation}
where $W$ is the unmixing matrix, $S$ are the reconstructed independent components, and $\lambda$ is a nongaussianity penalty parameter.
Gradients were computed using automatic differentiation with Zygote.jl~\cite{Innes:2019}.

\begin{figure}
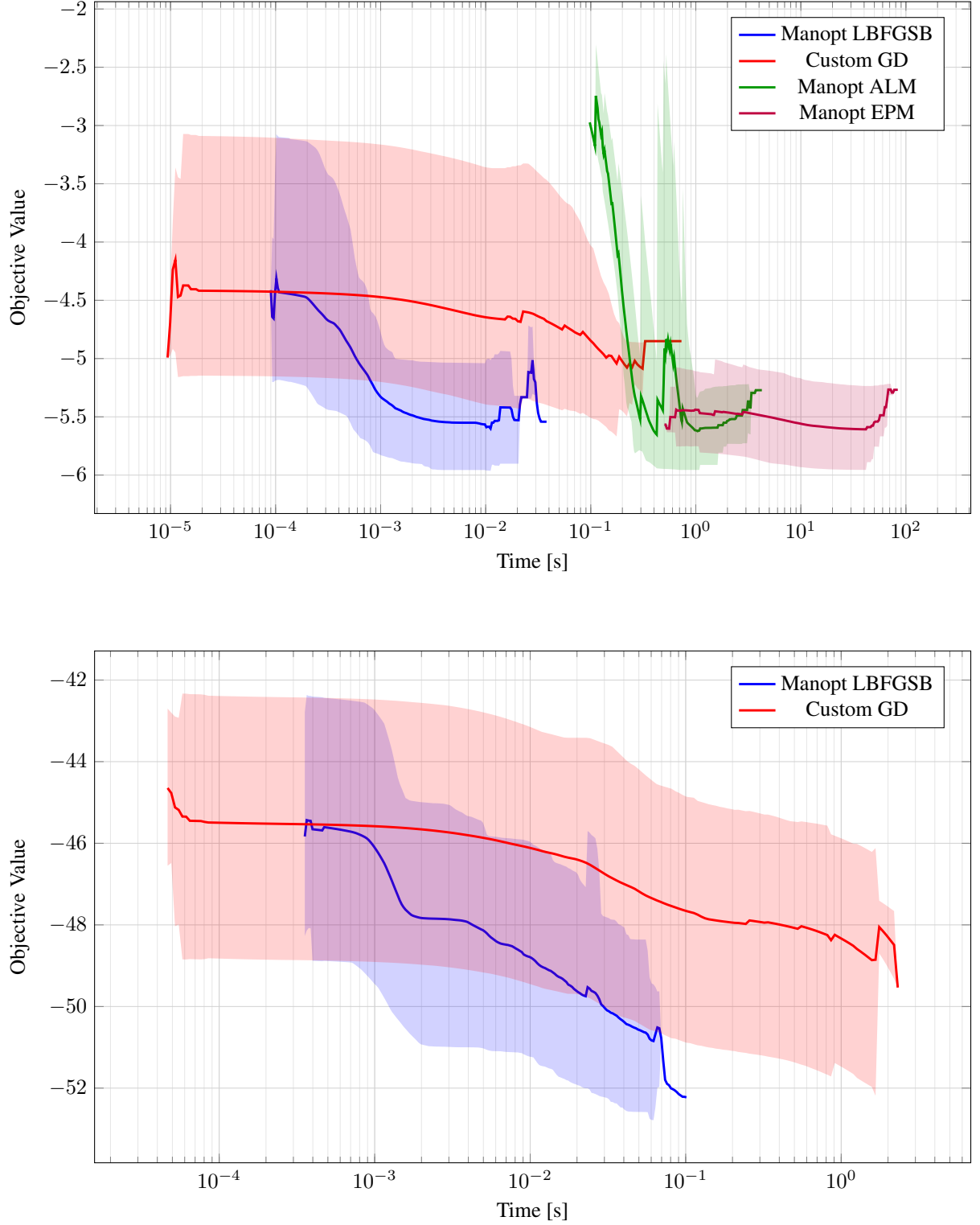

\begin{centering}

\include{objective_history_avg_ci_bss_ex1}

\include{objective_history_avg_ci_bss_ex2}

\end{centering}
\label{fig:bss_objective_history}
\caption{Objective history of different solvers for the amplitude-limited blind source separation problem (top plot: $n=50$, $r=3$, bottom plot: $n=500$, $r=3$). Series show the average over 20 random instances of the problem, and shaded areas represent 95\% confidence intervals.
The plots start after the first iteration, so the initial objective value is not shown.}
\end{figure}

Four solvers were considered: a simple gradient descent procedure with projection, Riemannian L-BFGS-B, Riemannian augmented Lagrangian method (ALM)~\cite{LiuBoumal:2020} and Riemannian exact penalty method (EPM)~\cite{LiuBoumal:2020}.
In the first experiments, we set $n=50$, $r=3$, $A=1$, and $\lambda=0.1$.
Averaged results over 20 random instances of the problem are shown in Figure~\ref{fig:bss_objective_history}.
Gradient descent performs iterations very quickly but it requires very small step size to avoid divergence.
ALM and EPM are more robust to the choice of parameters, but they require much more time per iteration and are only fast enough for very small problems.
Riemannian L-BFGS-B is the best performing solver overall.

\subsection{Bounded-variance maximum likelihood common principal components}

In the common principal components problem, we want to extend principal component analysis to data that belongs to $k$ classes with in such way that all classes share principal components.
Data from class $i$ is represented by a matrix $X_i \in \mathbb{R}^{r \times n_i}$ with $n_i$ samples and sample covariance matrix $S_i = \frac{1}{n_i} X_i X_i^\top$.
Additionally, we want to ensure that the variances are bounded between $d_{\min}$ and $d_{\max}$.
This leads to the following optimization problem:
\begin{equation}
  \operatorname*{arg\ min}_{Q \in \mathrm{SO}(r), \, D_i \in D(r)} \;  \sum_{i=1}^{k}
n_i \left\{
\log \bigl[\det(D_i)\bigr] + \operatorname{trace}
\left[ D_i^{-1} \odot \left(Q^{\top} S_i Q \right) \right]
\right\},
\end{equation}
where $Q$ is the common principal component matrix, $D_i$ are diagonal matrices representing the variances of the independent components for class $i$, $D(r)$ is the set of all diagonal matrices in $\mathbb{R}^{r \times r}$ with diagonal entries between $d_{\min}$ and $d_{\max}$, and $\odot$ represents the element-wise (Hadamard) product.

We have selected 19 standard datasets from the RDatasets repository~\cite{ArelBundock:2025} of different sizes, that possess clearly identifiable class column and at least 2 numerical features.
We have used the same boundaries for variances for all datasets, with $d_{\min} = 0.1$ and $d_{\max} = 10$.
The memory length was set to 2 for Riemannian L-BFGS-B, and the stopping criterion was set to $10^{-6}$ for the norm of the projected negative Riemannian gradient or 1000 iterations, whichever is reached first.
Table~\ref{tab:cpc_results} shows the results for Riemannian L-BFGS-B and Riemannian ALM.
Riemannian L-BFGS-B outperforms Riemannian ALM, achieving faster times with no constraint violations.
In most cases Riemannian L-BFGS-B is much faster than Riemannian ALM and achieves similar objective values to ALM when both return a feasible solution.
EPM solver is not included in the results as it performed significantly worse than ALM in preliminary experiments.

\begin{sidewaystable}
\centering
\caption{Detailed results for the bounded variance common principal components example.
Best values for each problem are shown in bold.}
\label{tab:cpc_results}

\begin{tabular}{
	l
  S[table-format=2.0]
  S[table-format=4.0]
  S[table-format=2.0]
  S[table-format=2.4]
  S[table-format=2.4]
  S[table-format=1.3e-2,tight-spacing=true,scientific-notation = true,]
  S[table-format=1.3e-2,tight-spacing=true,scientific-notation = true,]
  S[table-format=1.0]
  S[table-format=1.3e-2,tight-spacing=true,scientific-notation = true,]
}
\toprule
{\bfseries dataset} & {\bfseries features} & {\bfseries samples} & {\bfseries classes} & \multicolumn{2}{c}{{\bfseries time} [s]} & \multicolumn{2}{c}{\bfseries final cost} & \multicolumn{2}{c}{\bfseries violation} \\
 & & & & {L-BFGS-B} & {ALM} & {L-BFGS-B} & {ALM} & {L-BFGS-B} & {ALM}\\
\cmidrule(rl){1-1}\cmidrule(rl){2-2}\cmidrule(rl){3-3}\cmidrule(rl){4-4}\cmidrule(lr){5-6}\cmidrule(lr){7-8}\cmidrule(l){9-10}
COUNT/lbw & 9 & 189 & 2 & \bfseries 0.197 & 38.4 & 3.89e+06 & \bfseries 7.64e+04 & \bfseries 0 & 2.26e+05 \\
COUNT/medpar & 9 & 1495 & 3 & \bfseries 0.186 & 5.97 & \bfseries 1.01e+07 & 2.23e+17 & \bfseries 0 & 9.58e+05 \\
COUNT/rwm5yr & 16 & 19609 & 4 & \bfseries 0.586 & 15.1 & \bfseries 7.64e+09 & 7.9e+11 & \bfseries 0 & 4.12e+06 \\
Ecdat/Tuna & 7 & 13705 & 5 & \bfseries 0.336 & 11.1 & 9.84e+08 & \bfseries 1.11e+08 & \bfseries 0 & 8.62e+05 \\
HSAUR/bladdercancer & 2 & 31 & 2 & \bfseries 0.0571 & 0.447 & \bfseries 292 & \bfseries 292 & \bfseries 0 & \bfseries 0 \\
ISLR/Caravan & 85 & 5822 & 2 & \bfseries 7.16 & 540 & \bfseries -4.5e+04 & 4.96e+16 & \bfseries 0 & 207 \\
\midrule
KMsurv/kidney & 2 & 119 & 2 & \bfseries 0.134 & 0.682 & \bfseries 850 & \bfseries 850 & \bfseries 0 & 0.00381 \\
MASS/Pima.te & 7 & 332 & 2 & \bfseries 0.132 & 18.9 & \bfseries 4.04e+04 & \bfseries 4.04e+04 & \bfseries 0 & 1.07e-09 \\
MASS/Pima.tr & 7 & 200 & 2 & \bfseries 0.23 & 19.6 & \bfseries 2.61e+04 & \bfseries 2.61e+04 & \bfseries 0 & 8.74e-10 \\
MASS/Pima.tr2 & 4 & 300 & 2 & \bfseries 0.159 & 3.89 & \bfseries 2.65e+04 & \bfseries 2.65e+04 & \bfseries 0 & 1.1e-09 \\
MASS/biopsy & 8 & 699 & 2 & \bfseries 0.229 & 4.78 & \bfseries 8.03e+03 & \bfseries 8.03e+03 & \bfseries 0 & 4.44e-10 \\
MASS/birthwt & 10 & 189 & 2 & \bfseries 0.329 & 40.9 & 3.92e+06 & \bfseries 1.9e+06 & \bfseries 0 & 2.28e+05 \\
\midrule
datasets/iris & 4 & 150 & 3 & \bfseries 0.189 & 4.02 & \bfseries -789 & \bfseries -789 & \bfseries 0 & 1.63e-13 \\
rpart/stagec & 3 & 146 & 2 & \bfseries 0.503 & 2.34 & \bfseries 1.24e+03 & \bfseries 1.24e+03 & \bfseries 0 & \bfseries 0 \\
survey/nhanes & 5 & 8591 & 4 & 1.35 & \bfseries 0.906 & 4.74e+11 & \bfseries 3.87e+10 & \bfseries 0 & 1.04e+09 \\
survival/cancer & 3 & 228 & 2 & \bfseries 0.143 & 2.26 & 9.9e+05 & \bfseries 8.52e+03 & \bfseries 0 & 4.63e+04 \\
survival/heart & 7 & 172 & 2 & \bfseries 0.128 & 16.3 & 1.72e+06 & \bfseries 4.51e+05 & \bfseries 0 & 2.08e+05 \\
survival/kidney & 6 & 76 & 4 & \bfseries 1.01 & 19.3 & \bfseries 1.33e+05 & 1.97e+05 & \bfseries 0 & 3.13e+04 \\
survival/ovarian & 5 & 26 & 2 & \bfseries 0.176 & 8.81 & 1.64e+05 & \bfseries 2.55e+03 & \bfseries 0 & 8.6e+04 \\
\bottomrule
\end{tabular}
\end{sidewaystable}

\section{Conclusions}
\label{sec:conclusions}

In this paper we have presented Riemannian L-BFGS-B, a generalization of the L-BFGS-B algorithm to optimization problems on Riemannian manifolds with Euclidean bound constraints.
The algorithm is implemented in the Manopt.jl package, and we have shown that it can be used to effectively solve a wide variety of nonlinear problems.
In particular, we have shown that Riemannian L-BFGS-B can be used to solve large scale problems with bound constraints on the Euclidean space with minimal performance loss compared to pure Euclidean Fortran L-BFGS-B.
We suspect that the main reason for lower performance is omitting the subspace minimization step in our implementation.
This procedure is present in the Fortran L-BFGS-B and can significantly improve the performance when most of the variables are at the boundary of the feasible set, however it is difficult to generalize to Riemannian manifolds efficiently.

Benchmarks on two non-Euclidean problems show that Riemannian L-BFGS-B can be used to solve large scale problems with bound constraints and manifolds, and it significantly outperforms other general-purpose solvers for these problems like the augmented Lagrangian method and the exact penalty method.
Moreover, our solver guarantees that the iterates are always feasible, which is not the case for ALM and EPM.

The proposed method provides an effective and generic tool for solving problems with both manifold and bounded Euclidean variables, enabling a new class of problems to be efficiently solvable without the need for purpose-built solvers.

\section{Acknowledgements}

RB would like to thank MB and the AGH University for the kind hospitality and discussions during
a sabbatical stay in autumn 2025, where the main part of this paper stems from.
This research project supported by the program ,,Excellence initiative -- research university'' for the AGH University under the application IDUB 15636 (action D11).
\bibliographystyle{unsrtnat}
\bibliography{references}

@article{HuangGallivanAbsil:2015,
	title = {A Broyden Class of Quasi-Newton Methods for Riemannian Optimization},
	volume = {25},
	issn = {1052-6234},
	doi = {10.1137/140955483},
	pages = {1660--1685},
	number = {3},
	journaltitle = {{SIAM} Journal on Optimization},
	shortjournal = {{SIAM} J. Optim.},
	author = {Huang, Wen and Gallivan, K. A. and Absil, P.-A.},
	date = {2015-01},
	note = {Publisher: Society for Industrial and Applied Mathematics},
}

@article{HuangAbsilGallivan:2018,
    AUTHOR       = {Huang, Wen and Absil, P.-A. and Gallivan, K. A.},
    DOI          = {10.1137/17M1127582},
    JOURNAL      = {SIAM Journal on Optimization},
    NUMBER       = {1},
    PAGES        = {470--495},
    TITLE        = {A Riemannian BFGS method without differentiated retraction for nonconvex optimization problems},
    VOLUME       = {28},
    YEAR         = {2018},
}

@article{ZhuByrdLuNocedal:1997,
	title = {Algorithm 778: {L}-{BFGS}-{B}: {Fortran} subroutines for large-scale bound-constrained optimization},
	volume = {23},
	issn = {0098-3500},
	doi = {10.1145/279232.279236},
	number = {4},
	journal = {ACM Trans. Math. Softw.},
	author = {Zhu, Ciyou and Byrd, Richard H. and Lu, Peihuang and Nocedal, Jorge},
	month = dec,
	year = {1997},
	pages = {550--560},
}

@article{ByrdNocedalSchnabel:1994,
	title = {Representations of quasi-{Newton} matrices and their use in limited memory methods},
	volume = {63},
	issn = {1436-4646},
	doi = {10.1007/BF01582063},
	number = {1},
	journal = {Mathematical Programming},
	author = {Byrd, Richard H. and Nocedal, Jorge and Schnabel, Robert B.},
	month = jan,
	year = {1994},
	pages = {129--156},
}

@article{ByrdLuNocedalZhu:1995,
	title = {A {Limited} {Memory} {Algorithm} for {Bound} {Constrained} {Optimization}},
	volume = {16},
	issn = {1064-8275},
	doi = {10.1137/0916069},
	number = {5},
	journal = {SIAM Journal on Scientific Computing},
	author = {Byrd, Richard H. and Lu, Peihuang and Nocedal, Jorge and Zhu, Ciyou},
	month = sep,
	year = {1995},
	note = {Publisher: Society for Industrial and Applied Mathematics},
	pages = {1190--1208},
}

@misc{GrattonToint:2024,
	title = {{S2MPJ} and {CUTEst} optimization problems for {Matlab}, {Python} and {Julia}},
	doi = {10.48550/arXiv.2407.07812},
	number = {{arXiv}:2407.07812},
	publisher = {{arXiv}},
	author = {Gratton, Serge and Toint, Philippe L.},
	date = {2024-07-10},
	eprinttype = {arxiv},
	eprint = {2407.07812 [math]},
}

@article{GouldOrbanToint:2015,
	title = {{CUTEst}: a Constrained and Unconstrained Testing Environment with safe threads for mathematical optimization},
	volume = {60},
	issn = {1573-2894},
	doi = {10.1007/s10589-014-9687-3},
	shorttitle = {{CUTEst}},
	pages = {545--557},
	number = {3},
	journaltitle = {Computational Optimization and Applications},
	shortjournal = {Comput Optim Appl},
	author = {Gould, Nicholas I. M. and Orban, Dominique and Toint, Philippe L.},
	date = {2015-04-01},
}

@misc{Joyce:2010,
	title = {On manifolds with corners},
	url = {http://arxiv.org/abs/0910.3518},
	doi = {10.48550/arXiv.0910.3518},
	number = {{arXiv}:0910.3518},
	publisher = {{arXiv}},
	author = {Joyce, Dominic},
	urldate = {2024-04-17},
	date = {2010-10-13},
	eprinttype = {arxiv},
	eprint = {0910.3518 [math]},
}

@article{LiuWangZhao:2014,
	title = {Sparse {Covariance} {Matrix} {Estimation} {With} {Eigenvalue} {Constraints}},
	volume = {23},
	issn = {1061-8600},
	doi = {10.1080/10618600.2013.782818},
	number = {2},
	journal = {Journal of Computational and Graphical Statistics},
	publisher = {Taylor \& Francis},
	author = {Liu, Han and Wang, Lie and Zhao, Tuo},
	month = apr,
	year = {2014},
	pages = {439--459},
}

@article{ChenKangJianYishiYunpengDonald:2018,
	title = {Estimating large covariance matrix with network topology for high-dimensional biomedical data},
	volume = {127},
	issn = {0167-9473},
	doi = {10.1016/j.csda.2018.05.008},
	journal = {Computational Statistics \& Data Analysis},
	author = {Chen, Shuo and Kang, Jian and Xing, Yishi and Zhao, Yunpeng and Milton, Donald K.},
	month = nov,
	year = {2018},
	pages = {82--95},
}

@article{KramerSchaferBoulesteix:2009,
	title = {Regularized estimation of large-scale gene association networks using graphical {Gaussian} models},
	volume = {10},
	issn = {1471-2105},
	doi = {10.1186/1471-2105-10-384},
	journal = {BMC bioinformatics},
	publisher = {BMC Bioinformatics},
	author = {Krämer, N and Schäfer, J and Boulesteix, Al},
	month = nov,
	year = {2009},
}

@inproceedings{TuragaVeeraraghavanChellappa:2008,
	title = {Statistical analysis on stiefel and grassmann manifolds with applications in computer vision},
	doi = {10.1109/CVPR.2008.4587733},
	booktitle = {26th {IEEE} {Conference} on {Computer} {Vision} and {Pattern} {Recognition}, {CVPR}},
	author = {Turaga, Pavan and Veeraraghavan, Ashok and Chellappa, Rama},
	year = {2008},
	pages = {4587733},
}

@article{DaVeigaMarrel:2012,
	title = {Gaussian process modeling with inequality constraints},
	volume = {21},
	issn = {2258-7519},
	doi = {10.5802/afst.1344},
	number = {3},
	urldate = {2026-04-13},
	journal = {Annales de la Faculté des sciences de Toulouse : Mathématiques},
	author = {Da Veiga, Sébastien and Marrel, Amandine},
	year = {2012},
	pages = {529--555},
}

@article{SwilerGulianFrankelSaftaJakeman:2020,
	title = {A survey of {Constrained} {Gaussian} {Process} {Regression}: {Approaches} and {Implementation} {Challenges}},
	volume = {1},
	issn = {2689-3967, 2689-3975},
	shorttitle = {A {Survey} of {Constrained} {Gaussian} {Process} {Regression}},
	doi = {10.1615/JMachLearnModelComput.2020035155},
	language = {English},
	number = {2},
	journal = {Journal of Machine Learning for Modeling and Computing},
	publisher = {Begel House Inc.},
	author = {Swiler, Laura P. and Gulian, Mamikon and Frankel, Ari L. and Safta, Cosmin and Jakeman, John D.},
	year = {2020},
}

@article{PensoneaultYangXiu:2020,
	title = {Nonnegativity-enforced {Gaussian} process regression},
	volume = {10},
	issn = {2095-0349},
	doi = {10.1016/j.taml.2020.01.036},
	number = {3},
	journal = {Theoretical and Applied Mechanics Letters},
	author = {Pensoneault, Andrew and Yang, Xiu and Zhu, Xueyu},
	month = mar,
	year = {2020},
	pages = {182--187},
}

@article{WachterBiegler:2006,
	title = {On the implementation of an interior-point filter line-search algorithm for large-scale nonlinear programming},
	volume = {106},
	issn = {1436-4646},
	doi = {10.1007/s10107-004-0559-y},
	number = {1},
	journal = {Mathematical Programming},
	author = {Wächter, Andreas and Biegler, Lorenz T.},
	month = mar,
	year = {2006},
	pages = {25--57},
}

@unpublished{VanaretLeyffer:2026,
  author = {Vanaret, Charlie and Leyffer, Sven},
  title = {Implementing a unified solver for nonlinearly constrained optimization},
  year = {2026},
  note = {Accepted to Mathematical Programming Computation on Feb 22, 2026}
}

@BOOK{AbsilMahonySepulchre:2008,
  AUTHOR = {Absil, P.-A. and Mahony, R. and Sepulchre, R.},
  PUBLISHER = {Princeton University Press},
  DATE = {2008},
  DOI = {10.1515/9781400830244},
  TITLE = {Optimization Algorithms on Matrix Manifolds},
}

@BOOK{Boumal:2023,
  AUTHOR = {Boumal, Nicolas},
  PUBLISHER = {Cambridge University Press},
  URL = {https://www.nicolasboumal.net/book},
  DATE = {2023-03},
  DOI = {10.1017/9781009166164},
  TITLE = {An Introduction to Optimization on Smooth Manifolds},
}

@article{HyvarinenOja:2000,
	title = {Independent component analysis: algorithms and applications},
	volume = {13},
	issn = {0893-6080},
	shorttitle = {Independent component analysis},
	doi = {10.1016/S0893-6080(00)00026-5},
	number = {4},
	journal = {Neural Networks},
	author = {Hyvärinen, A. and Oja, E.},
	month = jun,
	year = {2000},
	pages = {411--430},
}

@article{LiuBoumal:2020,
	title = {Simple Algorithms for Optimization on Riemannian Manifolds with Constraints},
	volume = {82},
	issn = {1432-0606},
	doi = {10.1007/s00245-019-09564-3},
	pages = {949--981},
	number = {3},
	journaltitle = {Applied Mathematics \& Optimization},
	shortjournal = {Appl Math Optim},
	author = {Liu, Changshuo and Boumal, Nicolas},
	date = {2020-12-01},
}

@Manual{ArelBundock:2025,
  title = {Rdatasets: A collection of datasets originally distributed in various R packages},
  author = {Vincent Arel-Bundock},
  year = {2025},
  note = {R package version 1.0.0},
  url = {https://vincentarelbundock.github.io/Rdatasets},
}

@article{Innes:2019,
	title = {Don't {Unroll} {Adjoint}: {Differentiating} {SSA}-{Form} {Programs}},
	journal = {arXiv:1810.07951 [cs]},
	author = {Innes, Michael},
	month = mar,
	year = {2019},
	note = {arXiv: 1810.07951},
}

@article{LiuLeu:2025,
	title = {{ETCN}-{NNC}-{LB}: {Ensemble} {TCNs} {With} {L}-{BFGS}-{B} {Optimized} {No} {Negative} {Constraint}-{Based} {Forecasting} for {Network} {Traffic}},
	volume = {22},
	issn = {1932-4537},
	doi = {10.1109/TNSM.2025.3563978},
	number = {4},
	journal = {IEEE Transactions on Network and Service Management},
	author = {Liu, Jin-Xian and Leu, Jenq-Shiou},
	month = aug,
	year = {2025},
	pages = {3692--3704},
}

@article{FajemisinMaragnoDenHertog:2024,
	title = {Optimization with constraint learning: {A} framework and survey},
	volume = {314},
	issn = {0377-2217},
	doi = {10.1016/j.ejor.2023.04.041},
	number = {1},
	journal = {European Journal of Operational Research},
	author = {Fajemisin, Adejuyigbe O. and Maragno, Donato and den Hertog, Dick},
	month = apr,
	year = {2024},
	pages = {1--14},
}

@article{LaiYoshise:2024,
  title = {Riemannian Interior Point Methods for Constrained Optimization on Manifolds},
  volume = {201},
  doi = {10.1007/s10957-024-02403-8},
  number = {1},
  journal = {Journal of Optimization Theory and Applications},
  author = {Lai, Zhijian and Yoshise, Akiko},
  year = {2024},
  pages = {433–469}
}

@article{BergmannTenbrinck:2018,
	title = {A {Graph} {Framework} for {Manifold}-{Valued} {Data}},
	volume = {11},
	doi = {10.1137/17M1118567},
	number = {1},
	journal = {SIAM Journal on Imaging Sciences},
	author = {Bergmann, Ronny and Tenbrinck, Daniel},
	month = jan,
	year = {2018},
	pages = {325--360},
}

@misc{JuKoblerCollasKawanabeGuan:2026,
	title = {{SPD} {Matrix} {Learning} for {Neuroimaging} {Analysis}: {Perspectives}, {Methods}, and {Challenges}},
	doi = {10.48550/arXiv.2504.18882},
	publisher = {arXiv},
	author = {Ju, Ce and Kobler, Reinmar and Collas, Antoine and Kawanabe, Motoaki and Guan, Cuntai and Thirion, Bertrand},
	month = jan,
	year = {2026},
	note = {arXiv:2504.18882 [cs]},
}

@article{GongZhangYuanZhangXiong:2026,
	title = {An {Online} {Adaptive} {Physics}-{Constrained} {DMD} {Algorithm} with {Grassmann} {Manifold} {Spatial} {Mapping} for {Neutronic}-{Depletion} {Coupling} {Calculation}},
	issn = {0010-4655},
	doi = {10.1016/j.cpc.2026.110204},
	journal = {Computer Physics Communications},
	author = {Gong, Hanyuan and Zhang, Binhang and Yuan, Xianbao and Zhang, Yonghong and Xiong, Qingwen and Tang, Haibo and Zhou, Jianjun and Zhang, Sen and Xiao, Yunlong},
	month = apr,
	year = {2026},
	pages = {110204},
}

@article{FeiLiuJiaLiWeiChen:2025,
	title = {A {Survey} of {Geometric} {Optimization} for {Deep} {Learning}: {From} {Euclidean} {Space} to {Riemannian} {Manifold}},
	volume = {57},
	issn = {0360-0300},
	doi = {10.1145/3708498},
	number = {5},
	journal = {ACM Comput. Surv.},
	author = {Fei, Yanhong and Liu, Yingjie and Jia, Chentao and Li, Zhengyu and Wei, Xian and Chen, Mingsong},
	month = jan,
	year = {2025},
	pages = {123:1--123:37},
}

@article{MettesGhadimiKeller-ResselGuYeung:2024,
	title = {Hyperbolic {Deep} {Learning} in {Computer} {Vision}: {A} {Survey}},
	volume = {132},
	issn = {1573-1405},
	doi = {10.1007/s11263-024-02043-5},
	number = {9},
	urldate = {2026-05-06},
	journal = {International Journal of Computer Vision},
	author = {Mettes, Pascal and Ghadimi Atigh, Mina and Keller-Ressel, Martin and Gu, Jeffrey and Yeung, Serena},
	month = sep,
	year = {2024},
	pages = {3484--3508},
}

@article{ClosasOrtegaLesoupleDjuric:2024,
	title = {Emerging trends in signal processing and machine learning for positioning, navigation and timing information: special issue editorial},
	volume = {2024},
	issn = {1687-6180},
	doi = {10.1186/s13634-024-01182-8},
	number = {1},
	urldate = {2026-05-06},
	journal = {EURASIP Journal on Advances in Signal Processing},
	author = {Closas, Pau and Ortega, Lorenzo and Lesouple, Julien and Djurić, Petar M.},
	month = sep,
	year = {2024},
	pages = {84},
}

@inproceedings{TibermacineTibermacineZouaiRabehi:2024,
	title = {{EEG} {Classification} {Using} {Contrastive} {Learning} and {Riemannian} {Tangent} {Space} {Representations}},
	doi = {10.1109/ICTIS62692.2024.10894645},
	booktitle = {2024 {International} {Conference} on {Telecommunications} and {Intelligent} {Systems} ({ICTIS})},
	author = {Tibermacine, Ahmed and Tibermacine, Imad Eddine and Zouai, Meftah and Rabehi, Abdelaziz},
	month = dec,
	year = {2024},
	pages = {1--7},
}

@article{FrankeHefenbrockKoehlerHutter:2024,
	title = {Improving {Deep} {Learning} {Optimization} through {Constrained} {Parameter} {Regularization}},
	volume = {37},
	doi = {10.52202/079017-0286},
	language = {en},
	urldate = {2026-05-06},
	journal = {Advances in Neural Information Processing Systems},
	author = {Franke, Jörg K. and Hefenbrock, Michael and Koehler, Gregor and Hutter, Frank},
	month = dec,
	year = {2024},
	pages = {8984--9025},
}

@article{AxenBaranBergmannRzecki:2023,
  AUTHOR = {Axen, Seth D. and Baran, Mateusz and Bergmann, Ronny and Rzecki, Krzysztof},
  DOI = {10.1145/3618296},
  NUMBER = {4},
  TITLE = {{Manifolds.Jl: An Extensible Julia Framework for Data Analysis on Manifolds}},
  VOLUME = {49},
  YEAR = {2023},
  JOURNAL = {ACM Transactions on Mathematical Software},
}

@article{Bergmann:2022:1,
  AUTHOR = {Bergmann, Ronny},
  PUBLISHER = {The Open Journal},
  DOI = {10.21105/joss.03866},
  NUMBER = {70},
  PAGES = {3866},
  TITLE = {{Manopt.jl: Optimization on manifolds in Julia}},
  VOLUME = {7},
  YEAR = {2022},
  JOURNAL = {Journal of Open Source Software},
}

@article{LabsirLesoupleTourneret:2026,
	title = {K-{Means} and {Gaussian} {Mixture} {Models} on {Lie} {Groups}: {Application} to {Geometrical} {Clustering}},
	issn = {1941-0476},
	doi = {10.1109/TSP.2026.3678342},
	urldate = {2026-05-07},
	journal = {IEEE Transactions on Signal Processing},
	author = {Labsir, Samy and Lesouple, Julien and Tourneret, Jean-Yves},
	year = {2026},
	pages = {1--16},
}

@article{LezcanoIordachitaKim:2022,
	title = {Lie-{Group} {Theoretic} {Approach} to {Shape}-{Sensing} {Using} {FBG}-{Sensorized} {Needles} {Including} {Double}-{Layer} {Tissue} and {S}-{Shape} {Insertions}},
	issn = {1558-1748},
	doi = {10.1109/JSEN.2022.3212209},
	journal = {IEEE Sensors Journal},
	author = {Lezcano, Dimitri A. and Iordachita, Iulian I. and Kim, Jin Seob},
	year = {2022},
	pages = {1--1},
}

@article{BarTalmon:2024,
	title = {On {Interference}-{Rejection} {Using} {Riemannian} {Geometry} for {Direction} of {Arrival} {Estimation}},
	volume = {72},
	issn = {1941-0476},
	doi = {10.1109/TSP.2023.3322779},
	journal = {IEEE Transactions on Signal Processing},
	author = {Bar, Amitay and Talmon, Ronen},
	year = {2024},
	pages = {260--274},
}

@article{BouchardMalickCongedo:2018,
	title = {Riemannian {Optimization} and {Approximate} {Joint} {Diagonalization} for {Blind} {Source} {Separation}},
	volume = {66},
	issn = {1941-0476},
	doi = {10.1109/TSP.2018.2795539},
	number = {8},
	journal = {IEEE Transactions on Signal Processing},
	author = {Bouchard, Florent and Malick, Jérôme and Congedo, Marco},
	month = apr,
	year = {2018},
	pages = {2041--2054},
}

\end{document}